\documentclass[12pt]{amsart}
\usepackage{graphicx}
\usepackage{amssymb}
\usepackage{epstopdf}
\theoremstyle{plain}

 \theoremstyle{definition}
 
\numberwithin{equation}{section}

\title{DEGREE-$k$ LINEAR RECURSIONS Mod($p$) AND NUMBER FIELDS}
\author{T. MacHenry and Kieh Wong}

\begin{document}




\maketitle

\begin{quote} ABSTRACT:  Linear recursions of degree $k$ are determined by evaluating the sequence of Generalized Fibonacci Polynomials, $\{F_{k,n}(t_1,...,t_k)\}$  (isobaric reflects of the complete symmetric polynomials)  at the integer vectors $(t_1,...,t_k)$. If $F_{k,n}(t_1,...,t_k) = f_n$, then 
 $$f_n - \sum_{j=1}^k t_j f_{n-j} = 0,$$ 
 and $\{f_n\}$ is a linear recursion of degree $k$.
 On the one hand, the periodic properties of such sequences modulo a prime $p$ are discussed, and are shown to be related to the prime structure of certain algebraic number fields; for example, the arithmetic properties of the period are shown to characterize ramification of primes in an extension field. On the other hand, the structure of the semilocal rings associated with the number field is shown to be completely determined by Schur-hook polynomials.\\
 key words:  Symmetric polynomials, Schur polynomials, linear recursions, number fields.
 \end{quote}

\vspace{1.0cm}

{\large{\emph{1.INTRODUCTION}}}

A sequence $\{f_n\}$ is a \textit{linear recursion of degree $k$}, denoted by $[t_1, ...,t_k]$,
if, given a (finite) sequence of integers $t_1, ...,t_k$, the following equation is satisfied for all $n\in \mathbb{\mathrm{Z}}$: 

$$f_n - \sum_{j=1}^k t_j f_{n-j} = 0.$$

It is natural to ask about the periodic behaviour of linear recursions: given $k$ and $\mathbf{t} = (t_1, ...,t_k)$, is there an integer $c(\mathbf{t})$ such that $f_{n+c(\mathbf{t})} = f_n$ for all non-negative integers $n$? The answer is,  rarely  [8]. Theorem 2.1 below specifies the conditions under which a linear recursion \textit{is} periodic. 
A related question is,  when are such linear recursions periodic modulo a prime $p$. The answer in this case is, always.  The consequences and context of this answer is the content of Theorem 2.2 and the subject of this paper.     
The answers to these two questions are quite easy to  come by; but, it seems that their consequences are quite far-reaching and are implicated in two important areas of algebra:  the theory of symmetric polynomials, and the theory of  algebraic number fields.  It is the relation between these two areas that we wish to explore further.  The key to the connection lies in the fact that there is a bijective correspondence between the set of linear recurrences and the monic polynomials in   $ Z[X]$.  In particular,  there is a bijective correspondence between the linear recursions of degree $k$ and the monic polynomials of degree $k$,

 $$X^k -t_1 X^{k-1} -...-t_k$$ 
where

$${\bf t}=(t_1,...,t_k)$$

\noindent is an integer vector. We shall call this polynomial {\it the CORE polynomial}, and denote it by $ \mathcal{C}(t_1,...t_k) $. If we regard the components of the vector $\textbf{t}$ as indeterminates, then the core polynomial is the generic monic polynomial of degree $k$ with indeterminate coefficients, one for each degree $k$. 
$k$-th degree linear recursions are determined by a (linearly  recursive) sequence of polynomials in the variables $t_1,...,t_k$,  denoted by $\{F_n(\textbf{t})\}, \mathbf{t} = (t_1,...,t_k)$.  Explicitly, these polynomials satisfy the equation $$F_n(\mathbf{t}) - \sum_1^k t_jF_{n-j}(\mathbf{t)} =0.$$
\noindent This sequence of polynomials is a well known sequence which has been variously called the sequence of Generalized Fibonacci Polynomials (GFP or GFP(k)) \cite{TM},\cite{TM2},\cite{MT} and the $F-$ sequence of isobaric polynomials \cite{TM2}, \cite{MT}.  The sequence is just the sequence of complete symmetric polynomials in disguise \cite{MT}, \cite{Macd}.  The ring of symmetric polynomials in variables $\lambda_1,...,\lambda_k$ is isomorphic to the ring of isobaric polynomials \cite{TM2}, the isomorphism being given by $$t_j = (-1)^{j+1}\mathcal{E}_j,$$ 
where $\mathcal{E}_j$ is the $j-th$ Elementary Symmetric Polynomial in $k$ variables (ESP, or ESP(k)).  
This isomorphism clarifies the instrinsic connection between the core polynomials and properties of the linear recursions indexed by them.

In addition to being graded rings, these rings support several other algebraic structures; for example, they possess a convolution product, well known from multiplicative arithmetic functions \cite{TM}, and a Hopf algebra structure \cite{MT}.  Under the convolution product, the linearly recursive sequences form a (graded, commutative) group.  This product is closely connected with ordinary multiplication of the core polynomials.
Under the isomorphism given above, the power symmetric polynomials \cite{Macd} are mapped to the sequence of Generalized Lucas Polynomials (GLP or GLP(k)).  The justification for implicating Fibonacci and Lucas in this terminology is that  $F_{2,n}(1,1)$ and  $G_{2,n}(1,1)$ are just the usual Fibonacci and Lucas sequences of integers. When the GFP are mapped back to the ring of symmetric functions in the usual basis, they can be represented as

\begin{displaymath}\hat{F}_{k,n}(\lambda) = 
\frac{1}{\Delta} det \left( 
\begin{array}{c}
 1 \\ \lambda_j \\ ...\\ \lambda_j^{k-1}\\ \lambda_j^{k+n-1} \end{array} 
\right).
\end{displaymath}

\noindent$\Delta$ is the Vandermonde matrix, whose square is the discriminant of the core polynomial. For $k=2$ and $\textbf{t} = (1,1)$, these are the Binet formulae for the Fibonacci numbers \cite{TM2}.  The analogue of the Binet formulae for GLP(2) at $\textbf{t} = (1,1)$ is the sum of $n-th$ powers of the roots of the polynomial $X^2-X-1$. In general, the GLP, when written in terms of the roots of the core polynomial $\mathcal{C}(\textbf{t})$,  are sums of $n-th$ powers of the roots. 
Since the convolution product of two periodic ($mod(p)$) linearly recursive sequences is necessarily periodic ($mod(p)$), the period being the l.c.m. of the two periods, and since the convolution of two sequences is represented by a core which is a product of the two original core polynomials,  it suffices
to look at irreducible cores.  
 
 This leads us to number fields. Here we focus on two main points.  
 One is the particular way in which periodicity is reflected in the number fields.  The non-modular case---let us refer to this as $mod(1)$---involves cyclotomic fields.  In the $mod(p)$ case, the structures of interest are the semilocal rings derived from the number fields associated with the core polynomials.
The second point concerns the remarkable way in  which the well known classes of symmetric polynomials, the complete symmetric polynomials, the power symmetric polynomials,  the Schur polynomials, and, in particular, the weighted isobaric polynomials introduced in \cite{MT2}, are involved in the structure of the number fields and the semilocal rings associated with them (also see \cite{CL},\cite{HM},\cite{Lascoux2},\cite{Li},\cite{Prag}).  

The periodicity results for linear recursions mod(1) and mod($p$) are stated in Section 2 (Theorem 2.1, Theorem 2.2). However, their  proofs will be delayed until Section 3,  where they arise naturally as  consequences of the properties of the companion matrix of the core polynomial. Subsequently, we extend the results of Section 3 to the $mod(p)$ case in Section 4. In Section 5, we consider the extension fields of irreducible core polynomials over the rational numbers, their rings of integers,  and the resulting  $p$-semilocal rings. In the last section, we investigate the relation between the local prime number theory of the number fields and the $p$ -periodicity of the linear recursion implied by the core polynomial. In particular, we shall show that the prime number $p$ divides the $p$ -period of the recursion if and only if the prime $p$ ramifies in the extension field determined by the core polynomial. 

Finally, we stress the remarkable, intimate involvement of Schur polynomials in the structure and representation of the semilocal rings and in the number fields with which they are associated.

At the end of this paper,  we include a Maple Algorithm for computing the period of any $k$-order linear recursion modulo a prime $p$.  Professor Mike Zabrocki is the author of this algorithm.

{\large {\emph{2.PERIODIC LINEAR RECURSIONS}}}
 
 $\mathbf{THEOREM  \,2.1}$

A linear recursion is periodic if and only if every root of the core polynomial is a primitive complex root of unity.  In particular, if the core polynomial is the cyclotomic polynomial $CP(n)$ of degree  $\phi (n)$, where $\phi$ is the Euler totient function, then its associated linear recursion is periodic with period $n$ \cite{MT}.  (It is interesting to compare this theorem with the Lech-Mahler Theorem \cite{Cassels}.) $\square$
 
The proof will be discussed in Section 3.  

Denote the period of a linear recursion, either mod(p) or mod(1),  by  $c_p[t_1,...,t_k]$ where $p$ is either a rational prime or $p=1$, and the $t_j$ are the coefficients of the core polynomial.

$ \mathbf{THEOREM \,2.2}$

Every linear recursion is periodic modulo  $p$ for every rational prime $p$. The period  $c_p[\mathbf{t}] \leqslant p^k$. $\square$

This follows from simple combinatorial arguments, essentially the pigeonhole principle.  We improve this bound to a best bound, $p^k - 1$,  in Section 4. We observe that if a sequence is a periodic linear recursion, then
$F_{k,c_p} = F_{c_p} = 1, F_{c_{p-1}} = ...= F_{c_p -k+1} = 0,\, F_{c_p+1} = t_1.$

3.$\large{\emph{THE COMPANION MATRIX}}$

With each core polynomial, we associate its rational canonical matrix, the so-called \textit{companion matrix}.  We first consider the companion matrix for the generic core polynomial of degree $k$.

\begin{displaymath}
\mathbf{A} = 
\left(\begin{array}{cccc}
0&1&...&0\\
0&0&...&0\\
0&0&...&1\\
t_k&t_{k-1}&...&t_1
\end{array} \right) 
\end{displaymath}

Since $det{\bf A}= (-1)^{k+1} t_k$, $det{\bf A^n} = (-1)^{n(k+1)} t_k^n.  $ ${\bf A}$ is singular iff $t_k =0$. But if $t_k = 0$, the core polynomial is reducible;  so we assume  $\mathbf{A}$ to be non-singular.  
Thus $\mathbf{A}$ is invertible and generates a cyclic group ( finite, if the coefficients of the core polynomial satisfy the conditions of Theorem 2.1, otherwise, infinite).   The inverse of $\mathbf{A}$  is 

\begin{displaymath}
\mathbf{A^{-1}} = 
\left(\begin{array}{ccccc}
-t_{k-1}t_k^{-1}&-t_{k-2}t_k^{-1}&...&t_1t_k^{-1}&t_k^{-1}\\
1&0&...&0&0\\
0&1&...&0&0\\
...&...&...&...&...\\
0&0&...&1&0
\end{array} \right).
\end{displaymath}

Observe that the orbit of the k-th row vector of $\mathbf{A}$ under the action of $\mathbf{A}$ is just the first row of $\mathbf{A^2}$, and the action of $\mathbf{A^{-1}}$ on the first row of $\mathbf{A}$ is the second row of the inverse of $\mathbf{A}$.  So it is useful to consider the $ \infty \times k$ matrix whose row vectors are the elements of the doubly infinite orbit of  $\mathbf{A}$ acting on any one of them. 
For $k=3,$ $\mathbf{A^\infty}$ looks like this:
\begin{displaymath}
\mathbf{\mathbf{A^\infty}} =
\left(\begin{array}{ccc}
...&...&...\\
S_{(-n,1^2)} & -S_{(-n,1)} & S_{(-n)} \\
...&...&...\\
S_{(-3,1^2)} & -S_{(-3,1)} & S_{(-3)}\\
1&0&0\\
0&1&0\\
0&0&1\\
t_3&t_2&t_1\\
... & ... & ...\\
S_{(n-2,1^2)} & -S_{(n-2,1)} & S_{(n-2)}\\
S_{(n-1,1^2)} & -S_{(n-1,1)} & S_{(n-1)}\\
S_{(n,1^2)} & -S_{(n,1)} & S_{(n)}\\
... & ... & ...

\end{array}\right)_{\infty \times 3}
\end{displaymath}

\noindent This matrix has a number of important features which we summarize in

$\mathbf{THEOREM \, 3.1}$(cf.\cite{HM}, \cite{Lascoux}, \cite{Lascoux2}, \cite{CV1}, \cite{CV2})

(3.11) The row vectors  consist of the orbit of any row with $\mathbf{A}$ acting as a transformation matrix (on the right, say), and the components of the row vectors are just isobaric reflects of Schur-hook polynomials.

(3.12) The set  of $k\times k$ contiguous row vectors of $\mathbf{A^\infty}$,  with the entry in the lower right hand corner being $\mathbf{S_{(n)}}$, yields a (faithful) matrix representation of the cyclic group generated by $\mathbf{A}$:

  \begin{displaymath}
\mathbf{\mathbf{A^n}} = 
\left(\begin{array}{ccccc}
(-1)^{k-1}S_{(n-k+1,1^{k-1})}&...&(-1)^{k-j}S_{(n-k+1,1^{k-j})}&...&S_{(n-k+1)}\\
...&...&...&...&...\\
(-1)^{k-1}S_{(n,1^{k-1})}&...&(-1)^{k-j}S_{(n,1^{k-j})}&...&S_{(n).}
\end{array} \right)
\end{displaymath}

Or, more succinctly, we have  $${\bf A^n} = [(-1)^{k-j} \mathbf{S_{(i,1^{k-j})}}]_{k\times k},$$
where the entries are  isobaric Schur-hook reflects whose Young diagrams have arm length $i$ and leg length $k-j$.  

(3.13) The elements in each row of $\mathbf{A^\infty}$ are the coefficients of a representation of the powers (positive and negative) of any of the roots of the core polynomial---denoted by $\lambda^n$--- in terms of a basis consisting of the first $k-1$ powers of  $\lambda$:  $$\lambda^n = \sum_{j=0}^{k-1}(-1)^{k-j}\mathbf{S_{(n,1^{k-j})} }\lambda^{j}$$  for $n \in\mathbb{ Z}$,  where $\lambda$ is a root of the core polynomial, and the coefficients are Schur-hook reflects whose Young diagrams have  arm length $n$ and leg length $k-j$. 

(3.14) Each column of $\mathbf{A^\infty}$ is a $\textbf{t}$-linear recursion of Schur-hook polynomials. In particular, the right hand column is just the (doubly infinite) sequence of Generalized Fibonacci Polynomials, $\mathbf{F_{k,n}}$.

(3.15) $tr(\mathbf{A^n}) = \mathbf{G_{k,n}}(\textbf{t}$) for $n \in \mathbb{Z}$, where $\mathbf{G_{k,n}}$ is just the sequence of Generalized Lucas Polynomials, which is also a \textbf{t}-linear  recursion.

\underline{REMARK}  We note that the existence of the matrix $\mathbf{A^\infty}$ extends the sequences of Schur-hook polynomials, in particular, the GFP, as well as the GLP, in the negative direction.  It would be interesting to have a combinatorial interpretation of these negatively indexed symmetric functions.
One might compare this result with the theorem in \cite{TM2}, which gives rational convolution roots to all of the elements in the WIP-module \cite{MT2}, i.e., to all of the sequences of symmetric functions in the free $\mathbf{Z}$ -module generated by the Schur-hook polynomials. 

Proof (3.11-3.15).

(3.11) The orbit structure is a consequence of the construction of the matrix. Operation of the companion matrix on a $k$-vector of integers generates a linear recursion with respect to the vector $\mathbf{t}$. In fact,  the Schur-hook sequences sequences claimed in the theorem  \cite{MT}.

(3.12) follows from the arguments in (3.11).

(3.13)  follows from the Hamilton-Cayley Theorem.  A simple induction shows that these coefficients are just the stated Schur-hook functions of the theorem.

(3.14) This is discussed in (3.11). 

(3.15) The traces of the $k \times k$ -blocks are the sums of all of the Schur-hook (reflects) whose Young diagrams partition the same n;  but such sums of Schur-hooks are well known to be GLP of isobaric degree $n$  \cite{MT}. $\square$

The infinite companion matrix is a remarkable summary of all of the features connected with linear recursions (as enumerated in  Theorem 3.1):  It contains representations of the roots of the core polynomial as row vectors;  the right-hand column consists of GFP's, i.e., the generic k-th order linear recursions;  it displays the role of Schur-hook functions as both constituents of sequences of k-th order linear recursions---one of which is the GFP sequence--- and as  coefficients for a representation of the powers of the roots of the core polynomial. It contains a matrix representation of the free abelian group generated by the companion matrix, in particular,  a matrix representation of the free abelian group generated by any of the roots of the core.  It also contains, as traces, the GLP's.  Recalling that the GFP's and the GLP's are respectively, the isobaric versions of the complete symmetric polynomials and the power symmetric polynomials, we have connected the theory of linear recursion with two of the important bases of the algebra of symmetric polynomials. Moreover, we have introduced an extension of the symmetric polynomials to negatively indexed symmetric functions which are related to the reciprocals of powers of the roots of the core polynomial.  Thus we have a striking summary of the connection between the theory of equations and the theory of linear recursions within the ring of symmetric polynomials. 

$\mathbf{COROLLARY} \,3.2$ 

Given the k-th order linear recursion determined by the core $[t_1,\cdots,t_k]$,  with the companion matrix  $\mathbf{A}$, and denoting the cyclic group generated by $\mathbf{A}$ as $\mathbf{H}$,  we have that $\mathbf{H}$ is a finite cyclic group exactly when the linear recursion
is periodic, the order of the cyclic group $\mathbf{H}$ being the period of the recursion. Moreover, if the core polynomial is irreducible over the rationals, then every root of the core polynomial generates a finite cyclic group whose order is also the period of the linear recursion. 

Proof.   The proof follows immediately from Theorem 3.1. $\square$

Corollary 3.2 explains why Theorem 2.1 is true, for clearly, the only irreducible core polynomials having all of its roots periodic are those whose roots are roots of unity.

Applying the facts learned above about the companion matrix,  we now consider the periodic behaviour of linear recursions modulo a prime $p$. 
\vspace{0.50cm}

4.{\large{\emph{p-PERIODICITY AND THE COMPANION MATRIX}}}

Since the vector  $\mathbf{t} = [t_1,...,t_k]$ determines both the core polynomial and its associated linear recursion uniquely, we shall write  $[t_1,...,t_k]$ to denote either of these structures when the context is clear, and we shall extend the usage to the notation $[t_1,...,t_k]_p$ for a linear recursion $[t_1,...,t_k]$ modulo the prime $p$. As in Section 2, $c_p[t_1,...,t_k]$ denotes the period of $[t_1,...,t_k]_p$, and  $\mathbf{A_p}$ denotes the companion matrix with entries modulo $p$.  For any matrix $\mathbf{M}$, $tr\mathbf{M}$ denotes the trace of the $\mathbf{M}$.

THEOREM 4.1

(4.11)  $c_p[\mathbf{t} ] = c_p[t_1,...,t_k] \le p^k - 1. $
 
(4.12)  The (cyclic) group generated by $\mathbf{A_p}$ has order $c_p[\textbf{t}].$

(4.13) The columns of $\mathbf{A_p^\infty}$ have period $c_p[\textbf{t} ]$ .

(4.14) $\lambda^{c_p} [\textbf{t}] =_p 1$, where  $\lambda$ is a root of $\mathcal{C}{[\textbf{t}} ]$, and 
$c_p[\textbf{t]}$ is the least positive integer for which this is true;  i.e., $c_p[\textbf{t} ]$ is the $p$-order of $\lambda$.

(4.15) $tr\mathbf{A_p^n}$ is linearly recursive with period $c_p[\textbf{t}]$.

Proof (4.11-4.15).   (4.11) will be proved in  the next section.
Clearly $\mathbf{A_p}$ generates a cyclic group of order dividing $c_p[\mathbf{t}]$; on the other hand, since each of the columns of $\mathbf{A^\infty}$ is a linear recursion, they too must have a period $c_p[\mathbf{t}]$.

(4.14) is a direct consequence of (3.13),(3.14) and (4.13).

(4.15) is a consequence of (3.15). $\square$

\underline{REMARK}:  As pointed out above, Corollary 3.2 accounts for the truth of Theorem 2.1.  The core polynomials for the primitive $n-th$ roots of unity are the cyclotomic polynomials of degree $\phi(n)$,  whose roots have the obvious geometric period of $n$; that is,  $c_p[\textbf{t}]$ = $n$,  where $\textbf{t}$ is the appropriate vector of coefficients of the cyclotomic polynomial of degree $\phi(n)$.  This also affords a geometric interpretation of periodicity for the roots of the core polynomial in the plane of complex numbers with coordinates taken mod$(p)$,  which is analogous to the cyclotomic periodicity.
\vspace{1cm}

5.{ \large{\emph{THE NUMBER FIELD $\mathbf{R}[\mathbf{t}]$ \textsc{AND THE SEMILOCAL RING }$\mathbf{R_p}[\mathbf{t}]$}}}

$\mathbf{PROPOSITION} \,5.1$

 If the core polynomial  [\textbf{t}] is reducible, then $c_p[\textbf{t}]$ is the least common multiple of p-periods of its irreducible factors.

Proof.  This theorem follows easily from Proposition 3.13. $\square$

This suggests that we might as well consider only irreducible cores. But in that case,  we can also consider the number field $\mathcal{F}$ = $ \mathbb{Q}(\lambda) $ = $\mathbb{Q}[X]/id<\mathcal{C}(X)>$.
Let us denote the ring of integers (the maximal order) in this field by $\mathbf{R}[\textbf{t}]$ and we write $\mathbf{R} \otimes \mathbb{Z}_p = \mathbf{R_p}$.  

We can write the elements of the field $\mathcal{F}$ either as a module over the basis $\{1, \lambda,...,\lambda^{k-1}\}$, or uniquely as k-tuples $(m_o,...,m_{k-1})$ with entries from $ \mathbb{Q}$ with multiplication determined by the minimal polynomial of the field, or, as a result of the Hamilton-Cayley Theorem, as a module with the basis $\{\mathbf{I}, \mathbf{A},...,\mathbf{A^{k-1}}\}$. This gives a matrix representation of the elements in the field.  Call it the \textit{standard\textit{}} representation. We also have the same three options in $\mathbf{R_p}$ using these bases modulo $(p)$. Theorem (3.13) 
can be regarded as giving a representation of the powers of $\lambda$ in   $\mathcal{F}$, as polynomials in the integral $\lambda$ -basis where the coefficients are Schur-hook polynomials evaluated at $[\textbf{t}]$. Note that we have an induced \textit{standard} matrix representation in the ring $\mathbf{R_p}$.

One of the concerns of the theory of algebraic number fields is the relation between primes in the extension field $\mathcal{F}$ and the rational primes in $\mathbb{Z}$ that they sit over.  If we let $p$ be a rational prime generating the prime ideal $\mathbf{p}$ in  $\mathbb{Q}$, and let $\mathcal{P}$ be the ideal in $\mathbf{R}$ extending $\mathbf{p}$, then $\mathcal{P} = \mathcal{P}_{1}^{\epsilon_{1}}...\mathcal{P}_{s}^{\epsilon_{s}}$  is the prime decomposition of $\mathcal{P}$ in the Dedekind ring $\mathbf{R}$. If $f_j$ is the relative degree of the prime ideal $\mathcal{P}_{j}$, i.e., the degree of its minimal polynomial, then either $s = 1$ and  $\epsilon_1 = 1$, in which case $\mathcal{P}$ is a prime ideal, and $p$ is \textit{inert};  or, $s >1$ but $\epsilon_j = 1$ for all $j's$, in which case  $\mathcal{P}$ is the product of distinct prime ideals, and  $p$ \textit{splits};  or,  some $\epsilon_j > 1$ and $p$ \textit{ramifies}. These properties are reflected in the semilocal ring $\mathbf{R}_p$. Moreover,  there is a relation between the  phenomenon of periodicity of the linear recursion associated with the core polynomial and properties of the primes in the extensions of the core localized at $p$.  This will be discussed in the following sections. 

It is well known that for each core polynomial 
only a finite number of primes ramify;  when they do,  they divide the discriminant of the field. With few exceptions, the converse is also true, and those exceptions will not occur in our discussion \cite{J};  hence, for the purposes of this paper, $p$ ramifies if and only if $p|\Delta$,  where $\Delta$ is the discriminant of $\mathcal{F}$.  We shall want to use the following well known fact.

$\mathbf{PROPOSITION}\, 5.2$

$\Delta$ = $(-1)^{k(k-1)/2} \mathbf{ N(\mathcal{C}(X))}det\mathcal{C^{\prime}}(\mathbf{t})$.\;\;\quad\quad $ \square$

$\mathbf{\mathcal{C}'}(\mathbf{t})$ is the derivative of the core polynomial.

\noindent Noting that $ \mathbf{\mathcal{C'}}$ (the \textit{different}) can be regarded as an element of $\mathbf{R_p}$, and, denoting  $ \mathbf{\mathcal{C'}}$ by $\mathbf{D}$, we have

$\mathbf{COROLLARY}\, 5.3$

$ \mathbf{D_p} $ generates an ideal in $\mathbf{R_p}$ (the discriminant ideal) if and only if $p|\Delta$, that is, if and only if $p$ ramifies in $\mathbf{R}$. 

Proof.  $p$ divides the discriminant of the core polynomial  modulo $p$ if and only if $p$ ramifies,  which occurs if and only if the different vanishes modulo $p$ at a root of the core polynomial, and this happens if and only if the different generates an ideal in the semilocal ring $\mathbf{R_p}$ (the alternative being that the different is a unit in $\mathbf{R_p}$). $\square$

In keeping with the notation $\mathbf{A^\infty}[\mathbf{t}]$, we let $\mathbf{M^\infty}[\mathbf{t}]$ be the $\mathbf{H_p}$ -orbit of any row vector in the matrix $\mathbf{M}$.

Since, by construction, the columns of a standard matrix are $\mathbf{t}$-linear recursions, the following proposition can be proved by a simple induction.
\vspace{0.5cm}

$\mathbf{PROPOSITION}\, 5.4$

The right-hand column of $\mathbf{D^\infty}[\mathbf{t}]$ is the sequence of GLP's. $\square$

\noindent Here $\mathbf{D}$ is the standard matrix for the different, i.e., $\mathbf{D}= \mathbf{\mathcal{C'}}.$

6. {\large{\emph{STRUCTURE OF THE SEMILOCAL RING $\mathbf{R_p}$}}}

$\mathbf{R_p}$ is a finite, commutative ring;  it is, therefore, a semilocal ring whose structure is a well known part of classical algebraic number theory (e.g., \cite{J},\cite{McD}). 
We restate the structure theorem here (Theorem 6.4) for easy reference.  $\mathbf{R_p}$ also has an orbit structure under the action of the group generated by $\mathbf{A_p}$, which, while not mysterious, is not readily found in the literature, and plays an integral role in our results. We shall first discuss this orbit structure and then exploit the semilocal nature of $\mathbf{R_p}$.

If $t_k \neq 0\, mod(p)$, then $\mathbf{A_p}$ is non-singular, and, hence, is a unit in $\mathbf{R_p}$. The units in  $\mathbf{R_p}$ are exactly those elements with norms different from $0$, that is, having a standard matrix with non-zero determinant. An element with zero norm, then, either is zero or belongs to a proper  ideal.  Denote the group of units of $\mathbf{R_p}$ by $\mathbf{G_p}$ and its subgroup $gp<\mathbf{A_p}>$ by $\mathbf{H_p}$, the \textit{period subgroup}. Then  $\mathbf{R_p}$ is a $\mathbb{Z}_p (\mathbf{H_p})$-module, or more conveniently,  a right $\mathbf{H_p}$-module.  Clearly,  $\mathbf{R_p}$ is the disjoint union of its orbits under the action of $\mathbf{A_p}$. A number of observations follow from this fact. It will be useful to list them for future reference:

(1)  The orbit of zero is a singleton.

(2)  An ideal consists of the disjoint union of orbits, each of which has orbit length dividing $c_p[\mathbf{t}]$. (Clearly, two orbits are either disjoint or identical, up to cyclic permutation.)

(3) Two distinct orbits in the same maximal ideal differ from one another by a coset representative of  $\mathbf{H_p}$, i.e., if $O_1$ and $O_2$ are distinct orbits in the maximal ideal  $I$, then there is a coset representative $g$ of $\mathbf{H_p}$ in $\mathbf{G_p}$ such that $O_1g = O_2$. (Of course,  a coset representative may belong to the stabilizer of $\mathbf{H_p}$. $O_1$ and $O_2 $ need not be bijective.)

(4) The orbits of $\mathbf{G_p}$ are the cosets of $\mathbf{H_p}$.

(5) The columns of an orbit are $\mathbf{t}$-linearly recursive,  with a period dividing $c_p[\mathbf{t}]$. 

(6) The (standard) matrix representation of $\mathbf{R_p}$ is implicit in the orbit structure of $\mathbf{R_p}$.  If $\mathbf{m} \in \mathbf{R_p}$, and if $m_{i,j}$ is the $(i,j)-th$ component of the standard matrix representation $\mathbf{M}_p$ of  $\mathbf{m}$, the row vectors $\mathbf{m_i}$ of $\mathbf{M_p}$ are just the elements of the $\mathbf{A_p}$ -orbit  of $\mathbf{m}$. 

$\mathbf{PROPOSITION } \,6.1$

There is a one-to-one correspondence between  maximal ideals of $\mathbf{R_p}$ and  irreducible factors of $\mathcal{C }(X) mod(p) $.  $\square$

\vspace{0.50cm}

$\mathbf{PROPOSITION }\,6.2$ (Traces)

Let $\mathbf{m} = (m_0,...,m_{k-1}) \in \mathbf{R_p}$.  $tr(\mathbf{m}) = m_0 \mathbf{G_{k,0}} + ... +m_{k-1} \mathbf{G_{k,k-1}} $, where $\mathbf{\{G_{k,n}\}}$
is the sequence of Generalized Lucas Polynomials, i.e, the isobaric reflect of the complete symmetric polynomials. 

Proof. Express $\mathbf{m}$ as $(m_0,m_1,...,m_{k-1})$ and note that the rows of $\mathbf{M}$ are vectors $\mathbf{mA^i}$. Writing $\mathbf{A_j^i}$  for the $j-th$ column of $\mathbf{A^i}$, we have that the trace of $\mathbf{m}$ is
$$(\mathbf{m}A^0) \mathbf{A_1^0} + (\mathbf{m }\mathbf{A^1}) \mathbf{A_2^1} + ... +(\mathbf{m}\mathbf{A^{k-1}}) \mathbf{A_k^{k-1}}.$$
But a suitable rearrangement of this sum is just
$$m_0 tr\mathbf{A^0} + m_1 tr\mathbf{A^1} + ...+m_{k-1} tr\mathbf{A^{k-1}} .$$
which, by Theorem 3.15, yields Theorem 6.2. $\square$

Also note that since each component of a vector in an orbit is in exactly one trace computation, the sum of the components of vectors in an orbit is equal to the sum of the traces of the vectors in the orbit.  That is,

$\mathbf{PROPOSITION} \,  6.3$

(6.31) The sum of the elements of the $\mathbf{A_p}$ -orbit  of the vector $\mathbf{m}$ is the sum of the traces of the row vectors, $\mathbf{m_i}$, i.e.,
$$ \sum_{i,j} m_{i,j}= \sum_i tr \mathbf{m_i}.$$

(6.32)  If $\mathbf{m }\in \mathbf{R_p}$ , i.e, if  $det \mathbf{m}$ = $0$, then $$\sum_{orbits \, of\, I} \sum_i tr (\mathbf{m_i}) = 0. \qquad\qquad\square$$ 

$\mathbf{THEOREM} \,6.4$ (e.g.,\cite{McD}VI.2, \cite{J}). 

$\mathbf{R_p}$ is a semilocal ring. In particular, letting $\mathbf{J}(\mathbf{R_p}) = Rad(\mathbf{R_p}) = \mathbf{I_1} \cap ... \cap \mathbf{I_s }= \mathbf{I_1} ...\mathbf{I_s}$,  $\mathbf{I_1}, ...,\mathbf{I_s}$  a complete set of maximal ideals in $\mathbf{R_p}$,   there is a smallest integer  $m$ such that $\mathbf{J} = \mathbf{I_1}^m ...\mathbf{I_s}^m$, and $\mathbf{R_p}= \bigoplus_j \mathbf{R_p}/\mathbf{I_j^m}$,  where each factor is a local ring.  $\square$

REMARK.  For finite, commutative, semisimple rings, several of the radical operators coalesce.
The radical mentioned in the theorem can be taken, for example, to be the intersection of maximal ideals, or as the nilpotent radical.

We use the term 'p splits' to mean that the core polynomial factors modulo(p), but that it does not ramify.

$\mathbf{THEOREM} \,6.5$

(6.51) If $p$ is inert,  then $\mathbf{R_p}$ is a field.

(6.52) If $p$ splits, then $\mathbf{R_p}$ has a trivial radical, thus is semisimple, i.e., is the direct sum of  $s$ simple rings (fields), where $s$ is the number of prime ideals in the factorization of $\mathbf{R_p}$.

(6.53) If $p$ ramifies, then $\mathbf{R_p}$ has a non-trivial radical, and is a direct sum of  $s$ (non-trivial) local rings .  

PROOF.  Theorem 6.5 is a direct consequence of the structure theorem, Theorem 6.4. The $m$ in the theorem is the l.c.m. of the ramification indices.   $\square$

REMARK. An ideal element $ \mathbf{m}$ outside of the radical is cyclic, i,e, satisfies $\mathbf{m}^n = \mathbf{m}$ for some natural number $n$. If $\mathbf{m} = \mathbf{e}$ is an idempotent, then the powers of $\mathbf{eA}$ coincide with the orbit of $\mathbf{e}$. This is because $(\mathbf{eA})^n = \mathbf{eA}^n$;  thus $\mathbf{eA}$ generates a cyclic group of order dividing $c_p[\mathbf{t}]$.

Using the standard matrix representation of elements in $\mathbf{R}$ or in $\mathbf{R_p}$, we can assign to each element a \textit{rank} by letting  $rank\mathbf{m} = rank\mathbf{M}$, where $\mathbf{M}$ is the standard matrix representation of $\mathbf{m}$.  We then observe that all elements in the same orbit have the same rank;  that the rank of a unit is $k$, the degree of the core polynomial; and,  that the rank of an ideal element is at most the codegree of the ideal , i.e., $k-d$, where $d$ is the degree of the minimal polynomial of the ideal. (The rank of the representing matrix cannot exceed the degree of the minimal polynomial).

Denote the rank of an element $\mathbf{m}$ in $\mathbf{R_p}$ by $r(\mathbf{m})$.

$\mathbf{THEOREM \, 6.6}$

Suppose that $p$ splits and that $\{\mathbf{e_1},...,\mathbf{e_s}\}$ is a complete set of distinct primitive idempotents in $\mathbf{R_p}$. 
 $$r(\sum_1^s \mathbf{e_j}) =  \sum_1^s r(\mathbf{e_j}) = k.$$
 
Proof.  By (6.52),  $\mathbf{R_p}$ is semisimple.  We observe that:  $1 \le r(\mathbf{e_j}) < k$, and  since $\sum_1^s \mathbf{e_j} = 1$, $r(\sum_1^s \mathbf{e_j}) = k$.  The proof will then be a consequence of the following lemma and corollaries.

$ \mathbf{LEMMA}$ If we let $\mathbf{e}$ be the sum of the elements in any subset of the set of primitive idempotents $\{\mathbf{e_i}\}$, and let $\bar{\mathbf{e}}$ be the complementary sum, then $$r(\mathbf{e}) + r(\bar{\mathbf{e}})\leqslant k. $$ 
 
Proof.  If $\mathbf{E_1}$ and $\mathbf{E_2}$ are $ k \times k $ -matrices such that $\mathbf{E_1 E_2 }= \mathbf{0}$, then $r(\mathbf{E_1}) + r(\mathbf{E_2}) \leqslant k.$
 Since $\mathbf{E_1 E_2} = \mathbf{0}$ we have that $r(\mathbf{E_1}) \leqslant \nu(\mathbf{E_2})\leqslant k-r(\mathbf{E_2})$, where $\nu$ is the nullity of $\mathbf{E_2}$, and the lemma follows. $\square$

 $\mathbf{ COROLLARY \,6.61}$
\begin{center}
$r(\mathbf{e}) +  r(\mathbf{\bar{e}}) = k.$
\end{center}
Proof.  Using the above Lemma and the remark at the beginning of the proof of the theorem, we have $k = r(\mathbf{e} +  \mathbf{\bar{e}})  \leqslant   r(\mathbf{e}) +  r(\mathbf{\bar{e}}) = k.$  $\square$

 $\mathbf{ COROLLARY  \,6.62}$
\begin{center}
 $r(\mathbf{e_i} + \mathbf{e_j}) = r(\mathbf{e_i}) + r(\mathbf{e_j}).$
\end{center}
 
 Proof.  From Corollary 6.61, we have that  $r(\mathbf{e_1}) +  r(\mathbf{\bar{e_1}}) = k$, so that we can apply the above arguments to  $r(\mathbf{\bar{e_1}} ) = k -  r(\mathbf{e_1}) $ to deduce that $r(\sum_2^s\mathbf{e_i} ) = \sum_2^s r(\mathbf{e_i})$; hence, Corollary 6.62 holds.  $\square$
 
Theorem 6.6 follows now from the proof of Corollary 6.62. $\square$

$\mathbf{COROLLARY}$ 6.63  

\noindent   If we let $\mathbf{B_1}, ..., \mathbf{B_s}$ be the ideals  $\mathbf{R_p }\mathbf{\mathbf{e_1}}, ..., \mathbf{R_p} \mathbf{e_s}$ in $\mathbf{R_p}$, and let $\mathbf{B_j}^* = \mathbf{B_j }- \{\mathbf{0}\}$, then 
 $$\mathbf{B_1}^*\times...\times \mathbf{B_s}^* = \mathbf{\mathbf{G_p}},$$
 where $\mathbf{G_p}$ is the group of units of $\mathbf{R_p}$.
 If  $p$ does not ramify, the $\mathbf{B_j}^*$ are finite fields.

Proof.  This follows from Theorem 6.5, Theorem 6.6, the fact that ranks of non-zero elements of $\mathbf{R_p}$ are positive integers,  and that an element of $\mathbf{R_p}$ is a unit if and only if its norm is not zero  \cite{McD}. $\square$

$\mathbf{COROLLARY \,6.64}$
$$|\mathbf{G_p}| = |\mathbf{B_1}^*| ... |\mathbf{B_s}^*| = (p^{r_1}-1)...(p^{r_s}-1),$$ where $p^{r_i}$ is the order of $\mathbf{B_i}$ and $r_i$ is the rank of $\mathbf{e_i}$.     $\square$
  
$\mathbf{COROLLARY \,6.65}$

If $p$ splits,  the period $c_p[t_1,...,t_k] = lcm \{c_p(\mathbf{minpoly(e_i}))\}_1^s $. $\square$

The following result gives a remarkable connection between the p-periodicity of a linear recursion and the splitting properties of primes in associated rational number fields.
 
$\mathbf{THEOREM\,6.7}$   $p$ divides $c_p[\mathbf{t}]$ if and only if $p$ ramifies.

Proof.  First, we assume that $p | c_p[\mathbf{t}]$ and that $p$ does not ramify;  but then, by Theorem 6.5, $\mathbf{R_p}$ is semisimple, and, so by Corollary 6.64, $p | (p^{r_i}-1)$ for some $i$. A contradiction.  In particular, If $p$ does not ramify, $|\mathbf{G_p}|$ and $p$ are relatively prime.
In order to prove the converse,  we first prove the following lemma:

$\mathbf{LEMMA} $  If $\mathbf{e}$ is an idempotent in an ideal of $\mathbf{R_p}$, then the $\mathbf{H_p}$-orbit of $\mathbf{e}$ consists of the powers of $\mathbf{eA}$, a multiplicative cyclic group. In particular, the order of $\mathbf{eA}$ divides $c_p[\mathbf{t}]$.

Proof.  All of this follows easily from the fact that $ (\mathbf{e A)^{n}}   =  \mathbf{e}^{n} \mathbf{A^{n}}  = \mathbf{e A^n}$, that $\mathbf{eA^{c_p}} =\mathbf{ e}$,  and that the length of any orbit divides the period. $\square $

To finish the proof of the theorem, we observe that,  if $p$ ramifies, then each of the direct factors in $\mathbf{R_p}$ is a non-trivial local ring, say $\mathbf{B_j}$, where $\mathbf{B_j} =  \mathbf{I_j}/ {\mathbf{I_j}}^m$. If $\mathbf{B_j}^*$ is the group of units in $\mathbf{B_j}$,  then there is an idempotent  $\mathbf{e_j}$ in $\mathbf{I_j}$ and  $\mathbf{e_j }+ \mathbf{m}$ is a unit in the local ring $\mathbf{B_j}$, that is, is in $\mathbf{B_j}^{*}$,   whenever $\mathbf{m }\in {\mathbf{I}_\mathbf{j}^m}$, i.e. whenever $\mathbf{m }\in \mathbf{J}$. Moreover, there is a bijective correspondence between such $\mathbf{m}'s$ in $\mathbf{I_j}$ and the elements in the orbit of $\mathbf{e_j}$, so that  by the Lemma, $p$ divides $c_p$ .  Thus $p$ divides $|\mathbf{H_p}|$  and, hence, the order of $\mathbf{G_p}$. $\square$

 REMARK. The well known fact that a rational prime $p$ ramifies with respect to a cyclotomic extension over $CP(n)$  only if $p$ divides  $n$ now follows immediately from Theorems 2.1, 4.14, 6.7.

The following Maple algorithm for computing the period associated with a given core polynomial is due to Professor Mike Zabrocki of York University.
\begin{quote} \tt

\# companionmat - function which computes the companion matrix\\
\# to the polynomial $X^k - t_1 X^{k-1} - t_2 X^{k-2} - \cdots - t_k$\\
\# given the list [$t_1$, $t_2$, ..., $t_k$] of coefficients\\
> companionmat:=proc(tlist) local i,j;\\
    evalm([seq([seq(0,i=1..j),1, seq(0,i=1..nops(tlist)-1-j)],\\
     j=1..nops(tlist)-1), [seq(tlist[-i], i=1..nops(tlist))]]);\\
end:\\ \\
\# orbit- function which computes the orbit of an element \\
\# of $\bf R_p$ represented as a k-tuple, mlst,  when acted on \\
\# by the companion matrix with elements in ${\mathbb Z}_p$
corresponding\\
\#  to a list of coefficients, tlst\\
\# input: mlst - the sequence $(m_0, m_1, \ldots, m_{k-1})$\\
\#  tlst - the list of coefficients in the companion matrix\\
\#  p - a prime number\\
> orbit:=proc(mlst, tlst, p) local newlst,A,orb,domodm;\\
  domodm:=(E,m)->map(x->x mod m, E);\\
  A:=companionmat(tlst);\\
  newlst:=domodm(mlst,p);\\
  orb:=[evalm(newlst)];\\
  newlst:=domodm(evalm(newlst\&*A),p);\\
  while not convert(newlst,\`{}list\`{})=convert(orb[1],\`{}list\`{}) do \\
    orb:=[op(orb),evalm(newlst)];\\
   newlst:=domodm(evalm(newlst\&*A),p);\\
  od;\\
  orb;\\
end:\\ \\
\# return the size of the orbit computed in the last function\\
> lorbit:=proc(tlst, mlst, p);\\
  nops(orbit(tlst,mlst,p));\\
  end:\\
\end{quote}

\end{document}